\documentclass[12pt]{amsart}
\usepackage{latexsym, amsmath, amssymb, amsthm}
\usepackage{epsfig}
\usepackage{amscd}
\usepackage{latexsym}
\usepackage{pstricks,pst-node,cite}
\usepackage{eepic,epic}

\begin{document}
\baselineskip 16pt

\def \Aut {\mbox{\rm Aut\,}}
\def \Dic {\mbox{\rm Dic\,}}
\def \gcd {\mbox{\rm gcd\,}}
\def \Ker {\mbox{\rm Ker\,}}
\def \Cay {\mbox{\rm Cay\,}}
\def \Rot {\mbox{\rm Rot\,}}
\def \Mon {\mbox{\rm Mon\,}}
\def \PSL {\mbox{\rm PSL\,}}
\def \GL {\mbox{\rm GL\,}}
\def \AGL {\mbox{\rm AGL\,}}
\newfont{\sBbb}{msbm7 scaled\magstephalf}
\newcommand{\BC}{\mbox{\Bbb C}}
\newcommand{\BR}{\mbox{\Bbb R}}
\newcommand{\BZ}{\mbox{\Bbb Z}}
\newcommand{\sZ}{\mbox{\sBbb Z}}

\renewcommand{\labelenumi}{\theenumi)}
\newtheorem{theorem}{Theorem}[section]
\newtheorem{problem}{Problem}
\newtheorem{defin}{Definition}
\newtheorem{lemma}[theorem]{Lemma}
\newtheorem{prop}[theorem]{Proposition}
\newtheorem{conj}{Conjecture}
\newtheorem{op}{Open Problem}
\newtheorem{example}{Example}
\newtheorem{note}{Note}
\newtheorem{remark}{Remark}
\newtheorem{corollary}[theorem]{Corollary}
\newenvironment{pfA}{\medskip\noindent{\bf{Proof of Theorem 1.1:}}
  \hspace{-.4cm}      \enspace}{\hfill \qed \newline \smallskip}
  \newenvironment{pfB}{\medskip\noindent{\bf{Proof of Theorem 1.2:}}
  \hspace{-.4cm}      \enspace}{\hfill \qed \newline \smallskip}
  \newenvironment{pfC}{\medskip\noindent{\bf{Proof of Theorem 1.3:}}
  \hspace{-.4cm}      \enspace}{\hfill \qed \newline \smallskip}
\newcommand{\lra}{\longrightarrow}
\newcommand{\vect}[2]{\mbox{$({#1}_1,\ldots,{#1}_{#2})$}}
\newcommand{\comb}[2]{\mbox{$ \left( \begin{array}{c}
        {#1} \\ {#2} \end{array}\right)$}}
\setlength{\unitlength}{12pt}
\renewcommand{\labelenumi}{(\theenumi)}
\def\mod{\hbox{\rm mod }}

\author{Dongseok Kim}
\address{Department of Mathematics \\ Kyungpook National University \\ Taegu, 702-201 Korea}
\email{dongseok@knu.ac.kr}
\thanks{}

\author{Young Soo Kwon}
\address{Department of Mathematics \\Yeungnam University \\Kyongsan, 712-749, Korea}
\email{yskwon@yu.ac.kr}
\thanks{}

\author{Jaeun Lee}
\address{Department of Mathematics \\Yeungnam University \\Kyongsan, 712-749, Korea}
\email{julee@yu.ac.kr}
\thanks{This research was supported by the Yeungnam University research grants in 2007.}

\subjclass[2000]{05C10, 05C30}
\keywords{Cayley map, regular embedding}

\begin{abstract}
A Cayley map is a 2-cell embedding of a Cayley graph into an
orientable surface with the same local orientation induced by a
cyclic permutation of generators at each vertex. In this paper, we
provide classifications of prime-valent regular Cayley maps on
abelian groups, dihedral groups and dicyclic groups. Consequently,
we show that all prime-valent regular Cayley maps on dihedral
groups are balanced and all prime-valent regular Cayley maps on
abelian groups are either balanced or anti-balanced. Furthermore,
we prove that there is no prime-valent regular Cayley map
on any dicyclic group.
\end{abstract}

\title[prime-valent regular Cayley maps]{A classification of prime-valent regular Cayley maps on some groups}

\date{\today}
\maketitle

\section{ Introduction}

In this paper, we only consider undirected finite connected graphs
without loops or multiple edges.  For a simple graph $\Gamma$, an
arc of $\Gamma$ is an ordered pair $(x,y)$ of adjacent vertices of
$\Gamma$. Thus, every edge of $\Gamma$ gives rise to a pair of
opposite arcs. By $V(\Gamma)$, $E(\Gamma)$, $D(\Gamma)$ and $\Aut
(\Gamma)$, we denote the vertex set, the edge set, the arc set and
the automorphism group of $\Gamma$, respectively.  A graph
$\Gamma$ is said to be \emph{vertex-transitive},
\emph{edge-transitive} and \emph{arc-transitive} if $\Aut
(\Gamma)$ acts transitively on the vertex set, the edge set and
the arc set of $\Gamma$, respectively. A graph $\Gamma$ is
\emph{one-regular} if $\Aut (\Gamma)$ is arc-transitive and the
stabilizer of each arc in $\Aut(\Gamma)$ is trivial.  We consider
$\Aut(\Gamma)$ as an acting group on $V(\Gamma)$, $E(\Gamma)$ and
$D(\Gamma)$ according to the context.

For a simple graph $\Gamma$ with arc set $D$, an embedding of
$\Gamma$ or a \emph{map} with the underlying graph $\Gamma$ is  a
triple $\mathcal{M}=(D;R,L)$, where $R$ is a permutation of $D$
whose orbits coincide with the sets of arcs based at the same
vertex and $L$ is an involution of $D$ whose orbits are the pairs
of arcs induced by the same edge. The permutations $R$ and $L$ are
called a \emph{rotation} and an \emph{arc-reversing involution} of
$\mathcal{M}$, respectively. Let Mon$(\mathcal{M})$ be the
permutation group $\langle R,L\rangle$ generated by $R$ and $L$
and call it a \emph{monodromy group} of $\mathcal{M}$. Then,
Mon$(\mathcal{M})$ acts transitively on $D$.  Given two maps
$\mathcal{M}_1=(D_1;R_1,L_1)$ and $\mathcal{M}_2=(D_2;R_2,L_2)$, a
\emph{map isomorphism} $\phi : \mathcal{ M}_1 \rightarrow
\mathcal{ M}_2$ is a graph isomorphism between underlying graphs
such that $\phi R_1(x,y)=R_2\phi(x,y)$ for every arc $(x,y) \in
D_1$. In particular, if $\mathcal{ M}_1 = \mathcal{ M}_2 =
\mathcal{ M}$, then $\phi$ is called a \emph{map automorphism} on
$\mathcal{M}$. It follows that the automorphism group
$\Aut(\mathcal{M})$ of $\mathcal{M}$ acts semi-regularly on $D$.
If the action of $\Aut(\mathcal{M})$ on $D$ is regular then the
map itself is called \emph{regular}. It was shown that a map
$\mathcal{M}$ is regular if and only if the monodromy group of a
map $\mathcal{M}$ acts regularly on $D$~\cite{N}.

Let $G$ be a group and let $X=X^{-1}$ be a unit-free set of $G$
such that $\left< X \right> = G$. A \emph{Cayley graph}
$\Gamma=\Cay(G,X)$ is a graph with vertex set $G$ and two
vertices $g$ and $h$ are adjacent if and only if $g^{-1}h\in X$.
The set of left translations $\tilde G =\{L_g \mid g\in G \}$,
defined by $L_g(x)=gx$ forms a vertex-regular subgroup
of $\Aut(\Gamma)$.  A Cayley graph $\Cay(G,X)$ is \emph{normal} if the left regular
representation $\tilde{G}$ is a normal subgroup of $\Aut(\Cay(G,X))$.
Note that the arc set $D$ of the Cayley graph $\Cay(G,X)$ is
$\{(g, gx) \mid g\in G,\;x \in X \}$. Let $p$ be a cyclic
permutation of $X$. Then, a Cayley map $CM(G,X,q)$ is the map
$\mathcal{ M}=(D;R,L)$ with the rotation defined by
$R(g,gx)=(g,gq(x))$ and the arc-reversing involution $L$ defined
by $L(g, gx) = (gx, g)$ where $g \in G$ and $x \in X$. It is easy
to see that for every $g \in G$, $L_g R=R L_g$, hence $\tilde G$
is a subgroup of  $\Aut(\mathcal{M})$ acting regularly on
vertices. Furthermore, a Cayley map $\mathcal{ M}=CM(G,X,q)$ is
regular if and only if there exists an automorphism $\rho$ in the
stabilizer $(\Aut(\mathcal{M}))_v$ of a vertex $v$ cyclically
permuting the $|X|$ arcs based at $v$.  In this case,
$Aut(\mathcal{M})$ is a product of $\tilde G$ with a cyclic group
$\left< \rho \right> \cong \mathbb{Z}_n$, where $n=|X|$ (see
\cite{JA1,JAS}).

For an integer $t$, a Cayley map $\mathcal{ M} = CM(G,X,q)$ is \emph{
$t$-balanced} if $q(x)^{-1}=q^{t}(x^{-1})$ for every $x \in X$. In particular, a Cayley map $\mathcal{M}$ is
\emph{balanced} if it is $1$-balanced. A balanced Cayley map
$CM(G,X,q)$ is regular if and only if there exists a group
automorphism $\psi$ of $G$ whose restriction on $X$ is equal to
$q$~\cite{SS1}. In this case, the group $\tilde G$ is a normal subgroup of
$\Aut(\mathcal{ M})$ and $\Aut(\mathcal{ M})$ is a semidirect product of
$\tilde G$ by $\left< \psi \right>$. On the other hand, a Cayley map $\mathcal{M}$
is \emph{anti-balanced} if $\mathcal{M}$ is $(-1)$-balanced. For more general theory of Cayley
maps and their automorphisms, the reader is referred to
~\cite{JA1, JAS}.

For any positive integer $n$, $D_n = \left< a, b \mid a^n = b^2 = (ab)^2 = 1
\right>$ is the dihedral group of order $2n$ and $\Dic_n =\left<
a, b \mid a^{2n} = 1, \ b^2=a^n ~~\mbox{and}~~d^{-1}ab=a^{-1}
\right>$ is the dicyclic group of order $4n$. Their automorphism
groups are $$\Aut(D_n)= \{ \sigma_{i,j} \mid \sigma_{i,j}(a) =
a^i, ~ \sigma_{i,j}(b) =
 a^jb, ~~i,j \in \{1,2,\ldots,n \} ~~{\rm
and~~} (i,n) = 1 \},$$ and
$$\Aut(\Dic_n)= \{ \alpha_{i,j} \mid \alpha_{i,j}(a) =
a^i, ~ \alpha_{i,j}(b) =
 a^jb, ~~i,j \in \{1,2,\ldots,2n \} ~~{\rm
and~~} (i,2n) = 1 \}.$$
 In~\cite{CJT}, M. Conder et al.
developed a general theory of $t$-balanced Cayley maps and
classified regular anti-balanced Cayley maps on abelian groups.
In~\cite{KKF} and \cite{KO}, J.H. Kwak et al. classified regular
$t$-balanced Cayley maps on dihedral groups and dicyclic groups.
This paper focuses on the classification of prime-valent regular
Cayley maps on abelian groups, dihedral groups and dicyclic
groups. The main results are stated as the following three
theorems.

\begin{theorem} \label{main-abelian}
For any prime $p$, let $\mathcal{ M} = CM(G, X, q)$ be a
$p$-valent regular Cayley maps on abelian group $G$. Then
$\mathcal{M}$ is either balanced or anti-balanced. Moreover,
\begin{enumerate} \item[$(1)$] if $\mathcal{ M}$ is balanced, then
$\mathcal{ M}$ is isomorphic to $CM(\mathbb{Z}_{2}^{n}, X_1, q_1)$
for some elementary abelian $2$-group $\mathbb{Z}_{2}^{n}, X_1
\subset \mathbb{Z}_{2}^{n}$ and a cyclic permutation $q_1$ of $X_1$
such that
\begin{enumerate} \item[$(i)$] $X_1 = \{A^i \mathbf{x} \mid  0 \leq i \leq
p-1\; \}$ for some $A \in \GL_n(2)$ and $\mathbf{x} \in
\mathbb{Z}_{2}^{n}$ satisfying that $p$ is the smallest positive
integer such that $A^{p} = I$ and $\left< \mathbf{x}, A
\mathbf{x}, A^2 \mathbf{x}, \ldots, A^{p-1}\mathbf{x} \right> =
\mathbb{Z}_{2}^{n},$

 \item[$(ii)$]
 $q_1(A^i \mathbf{x}) = A^{i+1} \mathbf{x}$
for any $i \in \{0,1,\ldots,p-1\}$. \end{enumerate}

 \item[$(2)$]
If $\mathcal{ M}$ is anti-balanced, then $\mathcal{ M}$ is isomorphic
to
$$CM(\mathbb{Z}_{2p}, X_2 = \{ 1,\ 3, \ldots, \ 2p-1 \}, q_2 =
(1\;\;3\; \ldots \;2p-1)).$$
\end{enumerate}
\end{theorem}

\begin{theorem} \label{main2}
For any prime $p$, any $p$-valent regular Cayley map $\mathcal{ M} =
CM(D_n, X, q)$ on a dihedral group is balanced. Furthermore, $\mathcal{ M}$
is isomorphic to
$$CM(D_n, X_1=\left\{a^{\sum_{j=0}^t
\ell^j}b\mid 0\leq t\leq p-1\right\}, q_1 = (b \; ab \;
a^{\ell+1}b \; \ldots \; a^{\ell^{p-2}+\ell^{p-3}+\ldots+1}b))
$$ for some positive integer $\ell < n$ such that
$p$ is the smallest positive integer satisfying that $\ell^{p-1} +
\ell^{p-2} + \ldots + 1 \equiv 0 (\mod  n)$.
\end{theorem}

\begin{theorem} \label{maindic}
For any prime $p$, there is no $p$-valent regular Cayley map on a
dicyclic group.
\end{theorem}

Our paper is organized as
follows. In section~\ref{pre}, we review some known results on maps, Cayley maps and transitive permutation groups of prime order.
 In Section~\ref{main}, the classifications of prime-valent regular
 Cayley maps on abelian groups, dihedral groups and on dicyclic
groups are given. Furthermore, we give some remarks including that
for fixed prime $p$, the underlying Cayley graphs of $p$-valent
regular Cayley maps on dihedral groups  are one-regular with
finitely many exceptions.

\section{Preliminaries}\label{pre}

In this section, we review some results needed in the proofs  of our main theorems.

For a $k$-valent Cayley map $ \mathcal{M} = CM(\Gamma, X, q =
(x_1\;x_2 \ldots x_k) )$, let $\kappa$  be an involution on the
set $\{1,2, \ldots, k \}$ such that for any $i = 1,2, \ldots, k$,
$x_i^{-1} = x_{\kappa(i)}$ and call it the \emph{distribution of
inverses} of $\mathcal{M}$. We denote the group generated by two
permutations $(1\;2\; \ldots \; k)$ and $\kappa$ by $H(\kappa)$.

\begin{prop}[\cite{RSJTW}] \label{homo-fromCay-to-one}
Let $\mathcal{ M} = CM(G, X, q = (x_1\;x_2 \ldots x_k) ) =
(D;R,L)$ be a $k$-valent Cayley map  with the distribution of
inverses $\kappa$. Then, there exists a group epimorphism $f:
\Mon(\mathcal{ M}) \rightarrow H(\kappa)$ such that $f(R) =
(1\;2\; \ldots \; k)$, $f(L) = \kappa$ and $G$ is isomorphic to
$f^{-1}(H(\kappa)_{i})/\Mon(\mathcal{ M})_{e}$ for some arc $e \in
D$ and for some $i \in \{ 1,2,\ldots, k \}$ satisfying
$f(\Mon(\mathcal{ M})_{e}) \leq H(\kappa)_{i}$.
\end{prop}

Proposition~\ref{homo-fromCay-to-one} implies that if $\mathcal{
M} = CM(G, X, q = (x_1\;x_2 \ldots x_k) ) = (D;R,L)$ is a regular
$k$-valent Cayley map  with the distribution of inverses $\kappa$,
then there exists a group epimorphism $f: \Mon(\mathcal{ M})
\rightarrow H(\kappa)$ such that $G$ is isomorphic to
$f^{-1}(H(\kappa)_{i})$ for some $i, 1\le i\le k$. It means
that the group $H(\kappa)_i$ can be an epimorphic image of the
group $G$. Note that for any $i,j, 1\le i,j\le k$, two stabilizers
$H(\kappa)_i$ and $H(\kappa)_j$ are isomorphic because the group
$H(\kappa)$ acts transitively on the set $\{1,2,\ldots,k\}$.

The following proposition gives a classification of transitive
permutation groups of prime degree.

\begin{prop}[\cite{GNSS}] \label{prime-deg-pergp}
Let $G$ be a transitive permutation group of prime degree $p$.
Then, $G$ is isomorphic to one of the followings:
\begin{enumerate}
\item[$(i)$] the symmetric group $S_p$ or the alternating group
$A_p$;

 \item[$(ii)$] a subgroup of $\AGL_1(p)$,
1-dimensional
 affine group over the field $F$ of order $p$;

\item[$(iii)$] a permutation representation of $\PSL_2(11)$ of
degree $11$;

 \item[$(iv)$] one of the Mathiew groups
$M_{11}$ or $M_{23}$ of degree $11$ or $23$,
respectively;

 \item[$(v)$] a projective group $H$
with $\PSL_d(q) \leq H \leq {\rm P}\Gamma {\rm L}_d(q)$ of degree
$p = \frac{q^d-1}{q-1}$.
\end{enumerate}
\end{prop}

\begin{prop} [\cite{LG}] \label{projective-simple}
For any positive integer $d>1$ and for any prime power $q$, the
projective special linear group $\PSL_d(q)$ is simple except for
$\PSL_2(2)$ and $\PSL_2(3)$. Moreover,  $\PSL_2(2)
\cong S_3$, $\PSL_2(3) \cong A_4$ and $\PSL_2(4) \cong A_5$.
\end{prop}

M. \v{S}koviera and J. \v{S}ir\'{a}n
characterized regular balanced Cayley maps as follows.

\begin{theorem}[\cite{SS1}] \label{bal-characterize}
A Cayley map $\mathcal{ M} =CM(G, X, q)$ is a
regular balanced if and only if there exists an
automorphism $\rho$ of the group $G$ whose restriction on $X$ is
equal to $q$.
\end{theorem}

M. Conder et al. classified regular anti-balanced
Cayley maps on abelian groups. The classification implies the
following result.

\begin{prop}[\cite{CJT}] \label{abel-class-anti-regular}
For a fixed prime $p \geq 3$, let $\mathcal{ M} =CM(G, X, q)$ be a
$p$-valent regular anti-balanced Cayley map on abelian group.
Then, $\mathcal{ M}$ is isomorphic to the Cayley map
$CM(\mathbb{Z}_{2p}, X', q_1)$, where
 $$X' = \{ 1,\ 3, \ldots, \
2p-1 \}\quad \mbox{and} \;\; q_1 = (1\;\;3\; \ldots
\;2p-1).$$
\end{prop}

Note that the underlying graph in
Proposition~\ref{abel-class-anti-regular} is the complete
bipartite graph $K_{p,p}$.

In~\cite{KKF} and \cite{KO}, J.H. Kwak et al. classified regular
$t$-balanced Cayley maps on dihedral groups and dicyclic groups.
 The following
theorems find  classifications of regular anti-balanced Cayley
maps on dihedral groups and of regular $t$-balanced Cayley maps on
dicyclic groups.

\begin{theorem}[\cite{KKF}] \label{class-anti-regular}
Let $\mathcal{ M} =CM(D_n, X,q)$ be a regular anti-balanced Cayley map
on dihedral group  with $|X| \geq 3$. Then, $n$ is even number
$2n'$, $|X|=4$ and $\mathcal{ M}$ is isomorphic to a Cayley map
$CM(D_n, X=\{b,a,a^{2k}b,a^{-1}\},q=(b\; a \; a^{2k}b \; a^{-1}))$
for some $k$ satisfying $k^{2} \equiv -1 (\mod n')$.
\end{theorem}

\begin{theorem}[\cite{KO}] \label{class-t-regular-dic}
Let $\mathcal{ M} =CM(\Dic_n, X,q)$ be a regular $t$-balanced
Cayley map on dicyclic group  with $|X| \geq 3$. Then, $t=1$,
i.e., $\mathcal{ M}$  is balanced, and $\mathcal{ M}$ is
isomorphic to a Cayley map $CM(\Dic_n, X,q)$ with the cyclic
permutation $$q=(b\;\; ab\;\; a^{i+1}b \;\cdots \;
a^{i^{r-2}+i^{r-3}+\cdots+i+1}b\;\; a^nb \;\; a^{n+1}b \;\;
a^{n+i+1}b\; \cdots \;a^{n+i^{r-2}+i^{r-3}+\cdots+i+1})$$ on $X$,
where $1 \leq i \leq 2n-1$, $(i,2n)=1$, $r \geq 2$ and
$i^{r-1}+i^{r-2}+\cdots+i+1 \equiv n (\mod 2n)$.
\end{theorem}

For any positive integers $n, k$ and $\ell < n$ such that $k$ is
the smallest positive integer satisfying
$\ell^{k-1}+\ell^{k-2}+\ldots+1 \equiv 0 (\mod n)$, the Cayley map
$$CM(D_n, X=\{b, ab, a^{\ell+1}b,\ldots,
a^{\ell^{k-2}+\ell^{k-3}+\ldots+1}b \}, q=(b
\;\;ab\;\;a^{\ell+1}b\; \ldots
\;a^{\ell^{k-2}+\ell^{k-3}+\ldots+1}b))$$ is known to be a regular
balanced Cayley map~\cite{WY}. For a convenience, let
$\mathcal{T}$ be the set of triple $(n, \ell, k)$ of positive
integers such that $\ell < n$  and $k$ is the smallest integer
satisfying $\ell^{k-1} + \ell^{k-2} + \ldots + 1 \equiv 0 (\mod
n)$.  We denote the above Cayley map by $CM(n, \ell, k)$. Then,
the following theorem holds.

\begin{theorem}[\cite{WY}] \label{class-bal-regular}
Let  $\mathcal{ M} =CM(D_n, X, q)$ be a $k$-valent regular balanced
Cayley map. Then, the map $\mathcal{ M}$ is isomorphic to a Cayley map
$CM(n,\ell,k)$, namely,
\begin{eqnarray*}CM(D_n, X=\{b, ab, a^{\ell+1}b,\ldots,
a^{\ell^{k-2}+\ell^{k-3}+\ldots+1}b \}, q=(b
\;\;ab\;\;a^{\ell+1}b\; \ldots \;
a^{\ell^{k-2}+\ell^{k-3}+\ldots+1}b))\end{eqnarray*} for some
triple $(n, \ell, k) \in \mathcal{T}$. Moreover, for any two
triples $(n, \ell_1, k), (n, \ell_2, k) \in \mathcal{T}$, two
regular balanced Cayley maps $CM(n, \ell_1, k)$ and  $CM(n,
\ell_2, k)$ are isomorphic if and only if $\ell_1 = \ell_2$.
\end{theorem}

 Theorem~\ref{class-bal-regular} implies that for any fixed integer $k \geq 3$, $\{ CM(n, \ell,
 k) \mid (n, \ell, k) \in \mathcal{T} \}$ is the set of all
  $k$-valent regular balanced Cayley maps on the dihedral
 groups up to isomorphisms.

\section{Classifications of prime-valent Cayley maps}\label{main}

In this section, we classify prime-valent regular Cayley maps on
abelian groups,  dihedral groups and dicyclic groups,
respectively. First, we prove the following lemmas.

\begin{lemma} \label{affine}
For an odd prime $p$, let $G$ be a permutation group of degree $p$
 generated by $\rho=(1 \;\; 2 \;\;3\; \ldots \; p)$ and
 $\kappa$, where $\kappa(p) = p$ and $\kappa^2$ is the identity.
 If $G$ is isomorphic to a subgroup of $\AGL_1(p)$, then $\kappa$ is either the identity or $\kappa =
 (p)(1\;p-1)(2\;p-2)\ldots(\frac{p-1}{2} \; \frac{p+1}{2})$.
\end{lemma}
\begin{proof}
Because $\AGL_1(p)$ is a Frobenius group, $\left< \rho \right>$ is
a normal  subgroup of $G$.  Suppose that $\kappa(1)=i$. Then,
$\kappa \rho \kappa \in \left< \rho \right>$ and  $\kappa \rho
\kappa(p) = i$, which means $\kappa \rho \kappa = \rho^{i}$. It
implies that for any $k = 1, 2, \ldots, p$,
$$\kappa(k+i) =  \kappa \rho^{i}(k)  = \rho \kappa(k) =
\kappa(k)+1.$$ If $i= 1$ then $\kappa$ is the identity. Assume
that $i\neq 1$. Inserting $k=i$, one can get $\kappa(2i) = 2$.
And, taking $k = 2i$, one obtains $\kappa(3i) = 3$. By the same
process, we find $\kappa(ki) = k$ for any $k = 1, 2, \ldots, p$.
So, $\kappa(i^2) = i$. Since $\kappa(i)=1$ and $\kappa$ is an
involution,  $i^2 \equiv 1  (\mod p)$. Therefore, $i = p-1$ and
$\kappa =
 (p)(1\;p-1)(2\;p-2)\ldots(\frac{p-1}{2} \; \frac{p+1}{2})$.
\end{proof}

For an odd prime $p$ and a $p$-valent Cayley map $\mathcal{ M} =
CM(G, X, q)$ with the distribution of inverses $\kappa$, $\kappa$
fixes at least one element in $\{1,2, \ldots, p\}$. From now on,
we assume that $\kappa(p) = p$ without loss of generality.

 The next lemma shows that  any prime-valent regular Cayley
maps on abelian groups, dihedral groups or dicyclic groups  are either
balanced or anti-balanced.

\begin{lemma} \label{key-lemma}
For an odd prime $p \geq 3$, let $\mathcal{ M} = CM(G, X, p_1)$ be
a $p$-valent regular Cayley map  on an abelian group, a dihedral
group or a dicyclic group with the distribution of inverses
$\kappa$.
 Then, $\kappa$ is either the identity or $(p)(1\;p-1)(2\;p-2)\ldots(\frac{p-1}{2} \; \frac{p+1}{2})$.
\end{lemma}
\begin{proof}
By Proposition~\ref{homo-fromCay-to-one}, there exists a group
epimorphism $f: \Mon(\mathcal{ M}) \rightarrow H(\kappa)$ such
that $f(R) = (1\;2\; \ldots \; p)$, $f(L) = \kappa$ and $G$ is
isomorphic to $f^{-1}(H(\kappa)_{i})$ for some $i = 1,2,\ldots,
p$. It implies that the group $H(\kappa)_{i}$ is an epimorphic
image of $G$. Since $G$ is an abelian group, a dihedral group or a
dicyclic group, the group $H(\kappa)_{i}$ is either an abelian
group, a dihedral group or a dicyclic group. Moreover, since the
group $H(\kappa)$ generated by $(1\;2\; \ldots \; p)$ and $\kappa$
acts transitively on the set $\{1,2, \ldots, p\}$, $H(\kappa)$ is
isomorphic to one of the groups in
Proposition~\ref{prime-deg-pergp}.  \medskip

\emph{Case 1}:  $H(\kappa)$ is isomorphic to $S_p$ or $A_p$.

If $p=3$ and $H(\kappa) \cong S_3$ then $\kappa = (3)(1 \; 2)$. If
$p=3$ and $H(\kappa) \cong A_3$ then $\kappa$ is the identity.

Assume that $p > 3$ and $H(\kappa)$ is isomorphic to $S_p$ or
$A_p$. Then, for any $i =1, 2, \ldots, p$, the stabilizer
$H(\kappa)_{i}$ is isomorphic to  neither an abelian group, a
dihedral group
 nor a dicyclic group. Thus, it is impossible.\medskip

\emph{Case 2}:  $H(\kappa)$ is isomorphic to a subgroup of
$\AGL_1(p)$.

 By Lemma~\ref{affine}, $\kappa$ is either the identity or
 $(p)(1\;p-1)(2\;p-2)\ldots(\frac{p-1}{2} \; \frac{p+1}{2})$.
\medskip

\emph{Case 3}: $H(\kappa)$ is isomorphic to a permutation
representation of $\PSL_2(11)$ of degree $11$.

Because a permutation representation of $\PSL_2(11)$ of degree
$11$ is a transitive extension of a permutation representation of
$A_5$ of degree $10$, $H(\kappa)_{i}$ should be isomorphic to
$A_5$ for any $i = 1, 2, \ldots, 11$, which is isomorphic to
neither an abelian group, a dihedral group
 nor a dicyclic group.
\medskip

\emph{Case 4}: $H(\kappa)$ is isomorphic to the Mathiew groups
$M_{11}$ or $M_{23}$ of degree $11$ or $23$, respectively.

Since both the Mathiew groups $M_{11}$ of degree $11$ and $M_{23}$
of $23$ act 4-transitively on $11$-set and $23$-set, respectively,
 the stabilizer $H(\kappa)_{i}$ acts 3-transitively for any $i =1, 2, \ldots, p$.
 Thus, there exist $\alpha, \beta \in H(\kappa)_{i}$ such that $\alpha=(x \;\;y\;\;z) \cdots$
 and $\beta=(x \;\;y)(z) \cdots$. Note that the order of $\alpha$ is a multiple of $3$. If $H(\kappa)_{i}$
is isomorphic to  a dihedral group or a dicyclic group then
$\alpha$ should be a centralizer. But, one can easily check
$\alpha \beta \neq \beta \alpha$. Hence, $H(\kappa)_{i}$ can not be
isomorphic to an abelian group, a dihedral group nor a dicyclic
group.
\medskip

\emph{Case 5}: $H(\kappa)$ is isomorphic to a projective group $H$
with $\PSL_d(q) \leq H \leq {\rm P}\Gamma {\rm L}_d(q)$ of degree
$p = \frac{q^d-1}{q-1}$ for some prime power $q$.
\medskip

\emph{Subcase 5.1}: $d > 2$.  The stabilizer of the point $[0,0,
\ldots, 1]^{t}$ in $H$ contains a subgroup $S$ which is isomorphic
to $\PSL_{d-1}(q)$. Except for $(d, q) = (3,2)$ or $(3,3)$, the
group $\PSL_{d-1}(q)$ is simple by
Proposition~\ref{projective-simple}. For $d=3$ and $q=3$, the
group $\PSL_{2}(3)$ is isomorphic to the alternating group $A_4$.
For $d=3$ and $q=2$, let
\begin{gather*}
A = \begin{bmatrix} 1 & 1 & 0 \\ 0 & 1 & 0 \\ 1 & 1 & 1
\end{bmatrix} \quad \mbox{and} \quad B = \begin{bmatrix} 0 & 1 & 0 \\ 1 & 1 & 0 \\ 1 & 1 & 1
\end{bmatrix}.
\end{gather*}
Then, $A$ and $B$ are in the stabilizer of the point $[0,0,1]^{t}$
in $\PSL_3(2)$ and the orders of $A$ and $B$ are 4 and 3,
respectively. Note that an element of order $3$ in a dihedral
group or a dicyclic group is  a centralizer of the group. On the
other hand, one can easily check $AB \neq BA$. Hence, for any $i
=1, 2, \ldots, p$, the stabilizer $H(\kappa)_{i}$ can not be
isomorphic to an abelian group, a dihedral group nor a dicyclic
group. \smallskip

\emph{Subcase 5.2}: $d=2$ and $q$ is an odd prime power. Then, the
number $\frac{q^2 -1}{q-1} = q+1$ can not be prime. \smallskip

\emph{Subcase 5.3}: $d=2$ and $q$ is $2^{r}$ for some positive
integer $r$. If $q=2$ then $|X| = 3$. It means that $\kappa$ is
the identity or $(3)(1 \; 2)$. Now, we assume that $q$ is $2^{r}$
for some positive integer $r>1$. Then, there exist non-identity
elements $x, y \in GF(2^{r})$ whose orders are odd, where
$GF(2^{r})$ is the Galois field of order $2^{r}$. For these
elements $x, y \in GF(2^{r})$, let
\begin{gather*}
C = \begin{bmatrix}  1 & 0 \\ 0 & x
\end{bmatrix} \quad \mbox{and} \quad E = \begin{bmatrix} 1 & y \\ 0 &
x
\end{bmatrix}.
\end{gather*}
Then, $C$ and $E$ are in the stabilizer of the point $[1,0]^{t}$
in $\PSL_2(2^{r})$ and the order of $C$ is equal to that of $x$,
hence it is odd. Note that an element of odd order in a dihedral
group or a dicyclic group is a centralizer of the group. One can
easily check $CE \neq EC$. Thus, the stabilizer of the point
$[1,0]^{t}$ in $\PSL_2(2^{r})$ is isomorphic to neither an abelian
group,  a dihedral group nor a dicyclic group.
\medskip

Therefore, in any cases, $\kappa$ is either the identity or
$(p)(1\;p-1)(2\;p-2)\ldots(\frac{p-1}{2} \; \frac{p+1}{2})$.
\end{proof}

Note that if the distribution of inverses $\kappa$ of a Cayley map
$\mathcal{M} = CM(G, X, p_1)$ is the identity then $\mathcal{ M}$
is balanced and all elements (generators) in $X$ are involutions.
On the other hand, if $\kappa$ is
$(p)(1\;p-1)(2\;p-2)\ldots(\frac{p-1}{2} \; \frac{p+1}{2})$, then
the map $\mathcal{M}$ is anti-balanced.

\begin{pfA}
Let $\mathcal{ M} = CM(G, X, q)$ be a $p$-valent regular Cayley
map on an abelian group $G$ with the distribution of inverses
$\kappa$. Then, by Lemma~\ref{key-lemma}, $\kappa$ is either the
identity or $(p)(1\;p-1)(2\;p-2)\ldots(\frac{p-1}{2} \;
\frac{p+1}{2})$.

Assume that $\kappa$ is the identity. Then, $\mathcal{ M}$ is
balanced and all elements in $X$ are involutions. Furthermore, $X$
is an orbit under an automorphism $\psi$ of $G$ and the
restriction of $\psi$ on $X$ is $q$. Since $G$ is an abelian group
and $G$ is generated by involutions, $G$ is isomorphic to an
elementary abelian $2$-group $\mathbb{Z}_{2}^{n}$ for some
positive integer $n$. Note that the automorphism group of
$\mathbb{Z}_{2}^{n}$ is $\GL_n(2)$. Hence, $\mathcal{ M}$ is
isomorphic to $CM(\mathbb{Z}_{2}^{n}, X_1, q_1)$ such that
\begin{enumerate} \item[$(i)$] $ X_1 = \{A^i \mathbf{x}
\mid  0 \leq i \leq p-1\; \}$ for some $A \in \GL_n(2)$ and $
\mathbf{x} \in \mathbb{Z}_{2}^{n}$ satisfying that $p$ is the
smallest positive integer such that $A^{p} = I$ and $\left<
\mathbf{x}, A \mathbf{x}, A^2 \mathbf{x}, \ldots,
A^{p-1}\mathbf{x} \right> = \mathbb{Z}_{2}^{n}$ and
 \item[$(ii)$] $q_1(A^i
\mathbf{x}) = A^{i+1} \mathbf{x}$ for any $i \in
\{0,1,\ldots,p-1\}$. \end{enumerate}
 Next, we assume that $\kappa$ is
$(p)(1\;p-1)(2\;p-2)\ldots(\frac{p-1}{2} \; \frac{p+1}{2})$. Then,
$\mathcal{ M}$ is anti-balanced. Furthermore, by
Proposition~\ref{abel-class-anti-regular}, $\mathcal{ M}$ is
isomorphic to
$$CM(\mathbb{Z}_{2p}, X_2 = \{ 1,\ 3, \ldots, \ 2p-1 \}, q_2 =
(1\;\;3\; \ldots \;2p-1)).  $$
\end{pfA}

Note that the cardinality of the group $\GL_n(2)$ is
$$(2^n -1)(2^n - 2)(2^n - 2^2)\ldots(2^n -2^{n-1}) =
2^{1+2+\ldots+(n-1)}(2^n -1)(2^{n-1} - 1)(2^{n-2} - 1)\ldots(2
-1).$$ So, there exists an $A \in \GL_n(2)$
satisfying $A^p = I$ if and only if $p$ divide $2^k - 1$ for some
$1 \leq k \leq n$, namely, the order of $2$ in
the multiplicative group $\mathbb{Z}^*_p$ is at most $n$.

Next, we proceed to prove Theorems~\ref{main2} and \ref{maindic}.

\begin{pfB}
Let $\mathcal{M} = CM(D_n, X, q)$ be a $p$-valent regular Cayley
map on dihedral group with the distribution of inverses $\kappa$.
Then, by Lemma~\ref{key-lemma}, $\kappa$ is either the identity or
$(p)(1\;p-1)(2\;p-2)\ldots(\frac{p-1}{2} \; \frac{p+1}{2})$.

Suppose that $\kappa$ is $(p)(1\;p-1)(2\;p-2)\ldots(\frac{p-1}{2}
\; \frac{p+1}{2})$. Then, $\mathcal{ M}$ is anti-balanced. By
Theorem~\ref{class-anti-regular}, all regular anti-balanced Cayley
maps on dihedral groups are $4$-valent. So, it is impossible.
Therefore,  $\kappa$ is the identity and the map $\mathcal{ M}$ is
balanced. By Theorem~\ref{class-bal-regular}, it completes the
proof.
\end{pfB}

\begin{corollary} \label{prime-val-one-regular}
For any prime $p$, let $\Gamma = \Cay(D_n, X)$ be a $p$-valent
one-regular Cayley graph on dihedral group. Then, $\Gamma$ is
normal.
\end{corollary}
\begin{proof}
For the identity element $1$ in $D_n$, the stabilizer
$\Aut(\Gamma)_{1}$ is a cyclic group of order $p$ and acts
regularly on $N_1(1)$. Hence, there is a regular  Cayley map
$\mathcal{ M} = CM(D_n, X, q)$ whose underlying graph is $\Gamma =
\Cay(D_n, X)$. By Theorem \ref{main2}, $\mathcal{ M}$ is balanced.
Therefore, the left translation subgroup $\tilde{D_n}$ is a normal
subgroup of $\Aut(\Gamma) = \Aut(\mathcal{ M})$.
\end{proof}

\begin{pfC}
Let $\mathcal{M} = CM(\Dic_n, X, q)$ be a $p$-valent regular
Cayley map on a dicyclic group with the distribution of inverses
$\kappa$. Then, by Lemma~\ref{key-lemma}, $\mathcal{M}$ is either
balanced or anti-balanced.  By Theorem~\ref{class-t-regular-dic},
there is no regular anti-balanced Cayley maps on dicyclic groups
and every valancy of regular balanced Cayley maps on dicyclic
groups is even. Therefore, there is no prime-valent regular Cayley
map on dicyclic groups.
\end{pfC}

\begin{remark}
 For any prime $p$,
any $p$-valent one-regular Cayley graphs $\Gamma$ $= \Cay(D_n, X)$
on dihedral groups are underlying graphs of regular balanced
Cayley maps on dihedral groups by Theorem~\ref{main2} and
Corollary~\ref{prime-val-one-regular}. So, $\Gamma$ is isomorphic
to $\Cay(D_n,X$ $=\{b$, $ab$, $a^{\ell+1}b$, $\ldots$, $
a^{\ell^{p-2}+\ell^{p-3}+\ldots +1}b \})$ for some $\ell$ such
that $0 < \ell < n$ and $\ell^{p-1}$$+\ell^{p-2}$$+\ldots$$+1
\equiv 0 (\mod n)$. In~\cite{KKO}, it is shown that for
any prime $p$ and for any $n, \ell$ satisfying that $p$ is the
smallest positive integer such that $\ell^{p-1}$$+\ell^{p-2}$$ +
\ldots $$+ \ell $$+1 \equiv 0 $$(\mod n)$, the Cayley graph
$\Cay(D_n$,$X=\{b$, $ab$,$ a^{\ell+1}b$, $\ldots$,$
a^{\ell^{p-2}+\ell^{p-3}+ \ldots +1}b\})$ is one-regular except
finitely many such pairs $n$ and $\ell$. Namely, there exists a
constant $M$ which depends on $p$ such that for any $n > M$ and
$(n, \ell, p) \in \mathcal{T}$,  the Cayley graph
$\Cay(D_n$,$X=\{b, ab$, $a^{\ell+1}b$,$\ldots$,
$a^{\ell^{k-2}+\ell^{k-3} +\ldots +1}b\})$ is one-regular.  In~\cite{KKO}, it is also shown that for any $n \geq 31$ and
$(n, \ell, 5) \in \mathcal{T}$, the Cayley graph
$$\Cay(D_n,X=\{b, ab,
a^{\ell+1}b, a^{\ell^2 + \ell +1}b, a^{\ell^{3}+\ell^2+ \ell +
1}b\}) $$ is one-regular.
\end{remark}

By Theorems~\ref{class-bal-regular} and \ref{main2}, for a given
positive integer $n$ and for a given prime $p$, the number of
non-isomorphic $p$-valent regular Cayley maps on $D_n$ is the
number of positive integer $\ell$ satisfying $(n, \ell, p) \in
\mathcal{T}$. Note that if $\ell^{p-1} + \ell^{p-2} + \ldots + 1
\equiv 0 (\mod n)$ then $(\ell-1)(\ell^{p-1} + \ell^{p-2} + \ldots
+ 1) = \ell^{p} - 1 \equiv 0 (\mod n)$. So, the numbers $n$ and
$\ell$ are relatively prime.  For any positive integer $n =
2^ap_{1}^{a_{1}}p_{2}^{a_{2}}\ldots p_{t}^{a_{t}}~( p_{1}, p_{2},
\ldots, p_{t}$ are distinct odd prime numbers, $a \geq 0$ and
$a_{i}
> 0$ for each $i \geq 1$  $)$, the multiplicative group
$\mathbb{Z}_n^{*}$ is isomorphic to the product of  multiplicative
groups $\mathbb{Z}_{2^a}^{*} \times
\mathbb{Z}_{p_1^{a_1}}^{*}\times \ldots \times
\mathbb{Z}_{p_t^{a_t}}^{*}$ by Chinese remainder's Theorem.
Moreover, it is well known that $\mathbb{Z}_2^{*} = \{ 1 \}$,
$\mathbb{Z}_{4}^{*} \cong \mathbb{Z}_2$, $\mathbb{Z}_{2^a}^{*}
\cong \mathbb{Z}_2 \times \mathbb{Z}_{2^{a-2}}$ for $a \geq 3$ and
$\mathbb{Z}_{p_i^{a_i}}^{*} \cong \mathbb{Z}_{p_i^{a^i-1}(p_i-1)}$
for any odd prime $p_i$ with $a_i \geq 1$~\cite{JR}.

\begin{theorem} \label{main3}
Let $p$ be an odd prime and let  $n =
p_0^{a_0}p_{1}^{a_{1}}p_{2}^{a_{2}}\ldots p_{t}^{a_{t}}~( p_0,
p_{1}, p_{2}, \ldots, p_{t}$ are distinct primes and $p_0 = p$,
$a_0 \geq 0$ and $a_{i}
> 0$ for each $i \geq 1)$ be a positive integer. Then,
the number of non-isomorphic $p$-valent regular Cayley  maps on
the dihedral group $D_n$ is $(p-1)^{t}$ if $a_0=0$ or 1 and $p |
(p_i -1)$ for all $i =1,2, \ldots, t$; $0$ otherwise.
\end{theorem}
\begin{proof}
Let $\ell$ be a positive integer such that $\ell < n$ and
$\ell^{p-1} + \ell^{p-2} + \ldots + 1 \equiv 0 (\mod n)$. Then,
the numbers $n$ and $\ell$ are relatively prime and  $\ell^{p-1} +
\ell^{p-2} + \ldots + 1 \equiv 0 (\mod p_i^{a_i})$ for any $i=0,
1, \ldots, t$. Assume that  $\ell \equiv \ell_i (\mod p_i^{a^i})$
with $0 < \ell_i < p_i^{a_i}$ for any $i=0, 1, \ldots, t$. If
$\ell_i \equiv 1 (\mod p_i)$ for some $i=1, 2, \ldots, t$ then

$$\ell^{p-1} +
\ell^{p-2} + \ldots + 1 \equiv \ell_i^{p-1} + \ell_i^{p-2} +
\ldots + 1 \equiv p \neq 0 (\mod p_i), $$ which is a contradiction. Hence,
for any $i=1, 2, \ldots, t$, $\ell_i \neq 1 (\mod p_i)$. It
implies that the equation $\ell^{p-1} + \ell^{p-2} + \ldots + 1
\equiv 0 (\mod p_i^{a_i})$ is equivalent to the equation
$\ell_i^{p} \equiv 1 (\mod p_i^{a_i})$ for $i=1, 2, \ldots, t$. If
there exists an $i \in \{ 1,2, \ldots, t\}$ such that $p \nmid
(p_i -1)$ then there is no such $\ell_i$ because
$|\mathbb{Z}_{p_i^{a_i}}^{*}| = p_i^{a_i-1}(p_i-1)$. For any $p_i$
satisfying $p \mid (p_i -1)$, there exist $p-1$ elements $\ell_i$
such that $\ell_i \neq 1 (\mod p_i)$ and $\ell_i^{p} \equiv 1
(\mod p_i^{a_i})$.

Assume that  $a_0=1$. Then, for any $x \in \mathbb{Z}_p - \{1\}$,
$x^{p-1}+x^{p-2}+\ldots+ 1 = \frac{x^p-1}{x-1} \equiv 1 (\mod p)$.
Thus,  $\ell_0 = 1$ is the only integer such that $0 < \ell_0 < p$
and $\ell_0^{p-1} + \ell_0^{p-2} + \ldots + 1 \equiv  0 (\mod p)$.

Suppose that  $a_0 \geq 2$. For any  $x \in \mathbb{Z}_{p^{a_0}}$
with $x \neq 1 (\mod p)$, $x^{p-1}+x^{p-2}+\ldots+ 1 =
\frac{x^p-1}{x-1} \equiv 1 (\mod p)$. Hence, $\ell_0$ is $sp + 1$
for some $s$$(0 < s < p^{a_0-1})$. Moreover, $\ell_0$ is
$s'p^{a_0-1} + 1$ for some $s'$$(0 < s' < p)$ because $\ell_0$
should satisfy $\ell_0^{p} = (sp + 1)^p \equiv 1 (\mod p^{a_0})$.
But, for such a number $\ell_0$,
 \begin{eqnarray*} \ell_0^{p-1} + \ell_0^{p-2} & + &
\ldots + 1 \\ &\equiv& ((p-1)s'p^{a_0-1} + 1) + ((p-2)s'p^{a_0-1}
+
1) + \ldots + (s'p^{a_0-1} + 1) +1  \\
& \equiv & p (\mod p^{a_0}).    \end{eqnarray*}
Therefore, there exists no such a $\ell_0$.

Conversely, for any $\ell_0, \ell_1, \ldots, \ell_t$ such that  $0
< \ell_i < p_i^{a_i}$ and  $\ell_i^{p-1} + \ell_i^{p-2} + \ldots +
1 \equiv  0 (\mod p_i^{a_i})$ for any $i = 0,1,\ldots, d$, there
exists a unique $\ell$ such that $0 < \ell < n$, $\ell \equiv \ell_i
(\mod p_i^{a_i})$ and $\ell^{p-1}+\ell^{p-2}+\ldots+1 \equiv 0
(\mod n)$ by Chinese remainder's Theorem. Therefore, the number of
$\ell$'s satisfying  $0 < \ell < n$ and
$\ell^{p-1}+\ell^{p-2}+\ldots+1 \equiv 0 (\mod n)$ is $(p-1)^{t}$
if $a_0=0$ or 1 and $p | (p_i -1)$ for all $i =1,2, \ldots, t$; 0
otherwise.
\end{proof}

\end{document}